\documentclass[reqno]{amsart}[11pt,draft,amssymb]

\usepackage{amsthm,array,amssymb,amscd,amsfonts,latexsym}

\def\Q{{\mathbf Q}}
\def\Z{{\mathbf Z}}

\def\C{{\mathbf C}}
\def\R{{\mathbf R}}
\def\F{{\mathbf F}}
\def\H{{\mathbf H}}

\def\i{\mbox{\boldmath$i$}}
\def\j{\mbox{\boldmath$j$}}
\def\k{\mbox{\boldmath$k$}}
\def\1{\mbox{\boldmath$1$}}

\def\Tr{\mathrm{Tr}\hskip .2mm }
\def\tr{\mathrm{tr}}
\def\rk{\mathrm{rk}}
\def\End{\mathrm{End}}

\def\HDG{\mathrm{HDG}}

\def\I{\mathrm{id}}

\def\GL{\mathrm{GL}}
\def\SL{\mathrm{SL}}
\def\PSL{\mathrm{PSL}}

\def\M{\mathrm{M}}
\def\dim{\mathrm{dim}}
\def\O{{\mathcal O}}

\newtheorem{thm}{Theorem}[section]
\newtheorem{lem}[thm]{Lemma}
\newtheorem{cor}[thm]{Corollary}

\theoremstyle{definition}

\newtheorem{ex}[thm]{Example}

\newtheorem{rem}[thm]{Remark}

\title[Finite linear
groups, lattices, and products of elliptic curves]
{Finite linear groups, lattices, and\\ products of
elliptic curves}
\author[Vladimir  L. Popov]{Vladimir  L. Popov${}^*$}
\address{Steklov Mathematical Institute,
Russian Academy of Sciences, Gubkina 8, 119991 Moscow,
Russia} \email{popovvl\char`\@orc.ru}
\author[Yuri G. Zarhin]{Yuri G. Zarhin}
\address{Department of Mathematics, Pennsylvania State University,
University Park, PA 16802, USA
\newline \indent Institute for Mathematical Problems in
Biology, Russian Academy of Sciences, Pu\-shchino,
Moscow Region, Russia}
\email{zarhin\char`\@math.psu.edu}

\thanks{
 ${}^*$\,Supported by  Russian grants {
 \rus N{SH}--123.2003.01}, {\rus RFFI 05--01--00455}, and
Program of\\[-3pt] Mathematics Section of Russian Academy of
Sciences.}

\date{September 18, 2005}

\DeclareFontFamily{OT1}{wncyr}{\hyphenchar\font45 }
\DeclareFontShape{OT1}{wncyr}{m}{n}{%
  <5> <6> <7> <8> <9> gen * wncyr
   <10> <10.95> <12> <14.4> <17.28> <20.74>  <24.88>wncyr10}{}
\DeclareFontShape{OT1}{wncyr}{m}{it}{%
  <5> <6> <7> <8> <9> gen * wncyi
  <10> <10.95> <12> <14.4> <17.28> <20.74> <24.88> wncyi10}{}
\DeclareFontShape{OT1}{wncyr}{m}{sc}{%
  <5> <6> <7> <8> <9> <10> <10.95> <12> <14.4>
  <17.28> <20.74> <24.88>wncysc10}{}
\DeclareFontShape{OT1}{wncyr}{b}{n}{%
  <5> <6> <7> <8> <9> gen * wncyb
   <10> <10.95> <12> <14.4> <17.28> <20.74> <24.88>wncyb10}{}
\input cyracc.def
\def\rus{\usefont{OT1}{wncyr}{m}{n}\cyracc\fontsize{9}{11pt}\selectfont}

\DeclareMathSizes{9}{9}{7}{5} 

\begin{document}

\begin{abstract} Let $V$ be a finite
dimensional complex linear space and let $G$ be an
irreducible finite subgroup of $\GL(V)$.  For a
$G$-invariant lattice $\Lambda$ in $V$ of maximal rank,
we give a description of structure of the complex torus
$V/\Lambda$. In particular, we prove that for a wide
class of groups, $V/\Lambda$ is isogenous to a
self-product of an elliptic curve, and that in many
cases $V/\Lambda$ is isomorphic to a product of mutually
isogenous elliptic curves with complex multiplication.
We show that there are $G$ and $\Lambda$ such that the
complex torus $V/\Lambda$ is not an abelian variety, but
 one can always replace $\Lambda$ by another
 $G$-invariant lattice $\Delta$ such that $V/\Delta$ is
 a product of mutually isogenous
 elliptic curves with complex
 multiplication. We amplify these results with a
criterion, in terms of
 the character and the Schur $\Q$-index of
 $G$-module $V$, of
 the existence of a non\-ze\-ro $G$-invariant lattice in
~$V$.
\end{abstract}

\maketitle

\section{\bf Introduction}

This paper arose from the following observation made in
\cite{Popov}. Let $V$ be a complex linear space of
nonzero dimension $n<\infty$. Let $G\subset\GL(V)$ be a
finite irreducible {\it reflection group}, i.e., a
subgroup of $\GL(V)$ generated by (complex)
reflections, and let $\Lambda$ be a $G$-invariant
lattice in $V$ of rank $2n$ (hereinafter a {\it
lattice} is a discrete additive subgroup of a complex
or real linear space). All finite reflections groups
are classified in \cite{ST} (see also \cite{Cohen},
\cite{Popov}), and all lattices invariant with respect
to them are classified in \cite{Popov}.
 {\it A
posteriori}, it follows from this classification of
lattices that the complex torus $V/\Lambda$ is in fact
an abelian variety. Moreover, this classification
implies that $V/\Lambda$ is isogenous to a self-product
of an elliptic curve, and if $G$ is not
complexification of the Weyl group of an irreducible
root system, then $V/\Lambda$ is a product of mutually
isogenous elliptic curves with complex multiplication.

Our original goal was to give an {\it a priori},
independent of the classification of invariant lattices
proof of these properties of $V/\Lambda$. On this way,
we found out that in fact they hold for invariant
lattices of a much wider
 class of
irreducible finite subgroups of $\GL(V)$ than that of
reflection groups.  We give a description of
$V/\Lambda$ for every irreducible finite subgroup $G$
of $\GL(V)$ and $G$-invariant lattice $\Lambda$ in
$V$ of rank $2n$: the $G$-modules $V$ with Schur
${\bf Q}$-index~$1$ are considered in
Theorem~\ref{main2}, and the other $G$-modules $V$ in
Theorem~\ref{ratL}. In particular, we prove that in
the majority of cases (but not in all) $V/\Lambda$ is
an abelian variety; moreover, in many cases
$V/\Lambda$ is isogenous to a self-product of an
elliptic curve or even isomorphic to a product of
mutually isogenous elliptic curves with complex
multiplication. We show (Theorem~\ref{ratL} and
Example~\ref{non}) that $G$ and $\Lambda$ such that
the complex torus $V/\Lambda$ is not an abelian
variety do exist, but
 one can always replace $\Lambda$ by another
 $G$-invariant lattice $\Delta$ such that $V/\Delta$ is
 a product of mutually isogenous elliptic curves with complex
 multiplication (Theorem~\ref{deform}).
 We amplify these results with a criterion
 (in terms of
 the character and the Schur $\Q$-index of
 $G$-module $V$) of
 the existence of a non\-ze\-ro $G$-invariant
 lattice $\Lambda$ in
 $V$ (Theorem~\ref{criterion});
 the latter appears to be equivalent to the existence of
 $\Lambda$ of rank~$2n$.

Our approach hinges on the representation theory of
finite groups. However, for reflection groups, there is
another, geometric approach. It gives a key to the
classification of all invariant lattices and may be
useful for solving other problems (for instance, it
follows from the classification of invariant lattices
of reflection groups that elliptic curves arising from
$V/\Lambda$ for nonreal reflection groups
 are very specific:
 their endomorphism
 rings may be only the orders in
 $\Q(\sqrt{-d})$ for $d=1, 2,
 3,$ and $7$,
 see \cite{Popov};
 it would be
 interesting to find an {\it a
 priori}
 explanation of this phenomenon).
 Because of this reason,
 for reflection groups, we
 give
 in the last section
    a second, geometric proof
 of our result on the structure of $V/\Lambda$
 (Theorem~\ref{geom} and its proof).

\vskip 2mm

{\it Notation and terminology.}

$\Z$, $\Q$, $\R$, and $\C$  are respectively the ring
of integers, the field of rational numbers, the field
of real numbers, and the field of complex numbers.

$\H$ is the Hamiltonian quaternion $\R$-algebra $\bigl(\frac{-1,
-1}{\R}\bigr)$ (see, e.g.,\,\cite[\S 1.6]{PA}).

$\F_q$
is the finite field that consists of $q$ elements.

 The identity map of a set $S$ is
denoted by $\I_S$.

$\Z[S]$ is the subring of $\C$ generated by a subset
$S$ of $\C$.

For a subring $A$ of $\C$ and a subset $P$ of a linear
space $W$ over $\C$, the $A$-submodule of $W$ generated
by $P$ is denoted by $AP$.

$\Tr (P):=\{\tr\,(g)\mid g\in P\}$ where $P$ is a
subset of ${\rm End}_{\C}(V)$.

$\M_r(R)$ is the algebra of $r\times r$-matrices
 over a ring $R$ (associative and with identity element).

 $I_r$ is the identity matrix
 of $\M_r(R)$.


 $R^d$ is the space of column vectors over $R$ of height $d$.

 If $G$ is a finite subgroup
 of
$\GL(V)$, then $\chi{}^{}_{G, V}$ (respectively, ${\rm
Schur} {}^{}_{G, V}$) is the character (respectively,
the Schur index with respect to $\Q$) of $G$-module
$V$.

$Z_G$ is the center of $\Q$-algebra $\Q G$.

The field generated over $\Q$ by $\Tr (G)$ is denoted
by $\Q(\chi{}^{}_{G, V})$. If $\Q(\chi{}^{}_{G, V})$ is
$\Q$ (respectively, an imaginary quadratic number
field), then the character $\chi{}^{}_{G, V}$ is called
{\it rational} (respectively, {\it imaginary
quadratic}).

The $G$-module $V$ is called {\it orthogonal}
(respectively, {\it symplectic}) if there is a
symmetric (respectively, skew-symmetric) nondegenerate
$G$-invariant bilinear form $V\times V\to \C$.

\vskip 2mm

{\it Acknowledgement.} We thank E.~B.~Vinberg for
stimulating remarks. We learned
Corollaries~\ref{vbg} and \ref{vbgT} from him.


\section{\bf Some generalities}

Let $G$ be an irreducible finite subgroup of $\GL(V)$.
By Burnside's theorem, the irreducibility of $G$ is
equivalent to the equality
\begin{equation}\label{irr}
\C G=\End_{\C}(V).
\end{equation}

Since $G$ is a finite group, $\Q G$ is a
finite-dimensional $\Q$-algebra. The existence of
natural epimorphism of the group algebra of $G$ over
$\Q$ to $\Q G$ implies that $\Q G$ is a semisimple
$\Q$-algebra. Clearly, $Z_{G}$ is either a (number)
field or a product of (number) fields. Since $G\subset
\Q G$, the elements of $Z_{G}$ commute with $G$ and
hence, by \eqref{irr}, with $\End_{\C}(V)$. This
implies that
\begin{equation}\label{scal}
Z_{G} \subset \C\,\I_V.
\end{equation}
Hence $Z_{G}$ is a field. In turn, this implies that
$\Q G$ is a simple $\Q$-algebra and the\-re\-fore a
central simple $Z_{G}$-algebra. The latter means that
there is a central division $Z_{G}$-algebra $D$ and an
integer $r>0$ such that
\begin{equation}\label{D}
\Q G\simeq \M_r(D)\quad \mbox{(isomorphism of
$Z_{G}$-algebras)}.
\end{equation}

Below we shall naturally identify $\C\,\I_V$ with $\C$,
and $Z_{G}$ with the corresponding subfield of $\C$.
The above notation and conventions are kept throughout
the whole paper.

\begin{lem}\label{i}
\label{dim}  The natural $\C$-algebra homomorphism
\begin{equation}\label{psi}
\psi\colon \Q G\otimes_{Z_{G}}\!\C \longrightarrow \C
G=\End_{\C}(V)
\end{equation} is an isomorphism. In particular,
\begin{equation}\label{n2}
\dim^{}_{Z_{G}}(\Q G)= n^2.
\end{equation}
\end{lem}

\begin{proof} Since
 $\Q G\otimes_{Z_{G}}\!\C$ is a
 simple
$\C$-algebra, $\psi$ is injective. On the other hand,
\eqref{irr} implies that $\psi$ is surjective. \quad
$\square$
\renewcommand{\qed}{}\end{proof}

\begin{cor}\label{nr} $n^2=r^2\dim^{}_{Z_G}D$.
\end{cor}
\begin{proof} This follows from \eqref{D} and
\eqref{n2}. \quad $\square$
\renewcommand{\qed}{}\end{proof}

\begin{lem}
\label{Dtrace} $Z_{G}=\Q(\chi{}^{}_{G, V})$.
\end{lem}

\begin{proof} This is  well known
(see, e.g., \cite[Lemma 24.7]{D}). \quad $\square$
\renewcommand{\qed}{}\end{proof}

Recall  (see, e.g.,\,\cite[(70.4)]{CR}) that ${\rm
Schur} {}^{}_{G, V}=\min [K:\Q(\chi{}^{}_{G, V})]$ with
the mi\-ni\-mum taken over all subfields $K$ of $\C$
such that the linear group $G$ is defined over $K$. The
latter means that there exists a $G$-invariant $K$-form
of the $\C$-linear space $V$, i.e., a $K$-linear
subspace $L$ of $V$ such that $\dim_K(L)=n$ and $\C
L=V$. It is known (see, e.g.,\,\cite[(70.13)]{CR}) that
\begin{equation}\label{index}
{\rm Schur}^{}_{G, V}=\sqrt{\dim^{}_{\Q(\chi{}^{}_{G,
V})} D}.
\end{equation}
In particular, \eqref{index}, Corollary~\ref{nr}, and
Lemma~\ref{Dtrace} imply that
\begin{equation}\label{<=>}
\mbox{$Z_G$-algebras $\Q G$ and $\M_n(Z_G)$ are
isomorphic} \ \Longleftrightarrow \ {\rm
Schur}{}^{}_{G, V}=1.
\end{equation}

Clearly,  $\Z G$ is an order in $\Q G$; in particular,
$\Z G$ is a free $\Z$-module of rank $\dim_{\Q}(\Q G)$.
It is clear as well that $Z_{G}\cap \Z G$ is an order
in $Z_{G}$; in particular, it is a free $\Z$-module of
rank $\dim_{\Q}Z_G$.

\begin{lem}
{\rm(\cite[Section 3.1]{Popov})} \label{n2n} If there
exists a  nonzero $G$-invariant lattice $\Lambda$ in
$V$, then ${\rm rk}\,(\Lambda)=n$ or $2n$.
\end{lem}
\begin{proof}
Since $\Lambda$ is $G$-invariant, the $\C$-linear
subspaces $\C\Lambda=\R\Lambda+i\hskip .3mm\R\Lambda$
and $\R\Lambda\cap i\hskip .3mm\R\Lambda$ in $V$ are
$G$-invariant as well. The irreducibility of $G$ then
implies that $\C\Lambda=V$ and
$$
\R\Lambda\cap i\hskip .3mm\R\Lambda=\{0\}\hskip 2mm
\text{or}\hskip 2mm V.$$ Since
$\dim_{\R}\R\Lambda=\dim_{\R}i\hskip .3mm\R\Lambda={\rm
rk}(\Lambda)$, in the first case we obtain
$2n=\dim_{\R} V=\dim_{\R}\R\Lambda+ \dim_{\R}i\hskip
.3mm\R\Lambda=2{\rm rk}(\Lambda)$, so ${\rm
rk}(\Lambda)=n$. In the second case, $2n=\dim_{\R}
V=\dim_{\R}\R\Lambda={\rm rk}(\Lambda)$. \quad
$\square$
\renewcommand{\qed}{}\end{proof}

\begin{lem}
\label{nonrational0} If there exists a  nonzero
$G$-invariant lattice $\Lambda$ in $V$, then
$\chi{}^{}_{G, V}$ is either rational or imaginary
quadratic.
\end{lem}

\begin{proof} Pick a nonzero element $v\in
\Lambda$. Since $\Lambda$ is $G$-invariant, it is also
$\Z G$-invariant, and in particular, $Z_{G}\cap \Z
G$-invariant.  From \eqref{scal} we then deduce that
$(Z_{G}\cap \Z G) v\subseteq \Lambda\cap \C v$; whence
$Z_{G}\cap \Z G$ is a nonzero lattice in $\C$.
Therefore
 $\dim_\Q Z_{G}\leqslant 2$, i.e.,
$Z_{G}$ is either $\Q$ or a quadratic number field.
Since the orders of every real quadratic number field
 are not discrete in $\C$ (see, e.g.,
 \cite[Ch.\,II, \S 7]{BS}), the claim now follows from
 Lemma~\ref{Dtrace}.
\quad $\square$
\renewcommand{\qed}{}\end{proof}


\begin{lem}
\label{nonrational} Suppose that $\chi^{}_{G,V}$ is not
rational. If there exists a nonzero $G$-in\-va\-ri\-ant
lattice $\Lambda$ in $V$, then
\begin{itemize} \item[${\rm(i)}$]
$Z_{G}$ is an imaginary quadratic number field;
\item[${\rm(ii)}$] ${\rm Schur}{}^{}_{G, V}=1$;
\item[${\rm(iii)}$] $\rk(\Lambda)=2n$.
\end{itemize}
\end{lem}

\begin{proof}
Since $\chi{}^{}_{G, V}$ is not rational,
 ${\rm(i)}$ readily follows
 from
 Lemmas
 \ref{Dtrace} and \ref{nonrational0}.
The $\Q$-linear space $\Q\Lambda$ carries a natural
structure of $Z_G$-linear space, and (i) imples that
$\dim_{Z_G}(\Q\Lambda)=\dim_{\Q}(\Q\Lambda)/2=
\rk(\Lambda)/2$. Since $\Q\Lambda$ is $\Q G$-stable, we
get a $Z_G$-algebra homomorphism $\varphi:\Q G \to
\End_{Z_G}(\Q\Lambda)$. Since $\varphi(1)={\rm
id}_{\Q\Lambda}$ and $\Q G$ is simple, $\varphi$ is
injective. This and \eqref{n2} then imply
\begin{equation}\label{<<}
n^2=\dim_{Z_G}(\Q G) \leqslant
(\dim_{Z_G}(\Q\Lambda))^2= (\rk(\Lambda)/2)^2 \leqslant
n^2.
\end{equation}
From \eqref{<<} we deduce that ${\rm(iii)}$ holds and
$\varphi$ is an isomorphism. The latter property and
\eqref{<=>} clearly imply ${\rm(ii)}$. \quad $\square$
\renewcommand{\qed}{}\end{proof}

\begin{cor}
[E.~B.~Vinberg]
\label{vbg} Suppose that the
$G$-module $V$ is not self-dual. If there exists a
nonzero $G$-invariant lattice $\Lambda$ in $V$, then
the conclusions of Lemma~{\rm \ref{nonrational}} hold.
\end{cor}

\begin{proof}
Since $V$ is not self-dual, the character
$\chi{}^{}_{G, V}$ is not real valued and therefore is
not rational. The claim now follows from Lemma
\ref{nonrational}. \quad $\square$
\renewcommand{\qed}{}\end{proof}

\begin{lem}
\label{split} Suppose that the greatest common divisor
of the integers $\ker(u)$, where $u$ runs through $\Q
G$, is equal to $1$. Then ${\rm Schur}{}^{}_{G, V}=1$.
\end{lem}
\begin{proof}
Since $G$ is irreducible and the elements of $\Tr\,(G)$
are integral algebraic numbers, \cite[30.10]{CR}, the
assumption of lemma implies by \cite[Lemma 3]{Vinberg}
and Lemma~\ref{Dtrace} that $G$ is defined over $Z_G$,
whence the claim by \eqref{<=>}. \quad $\square$
\renewcommand{\qed}{}\end{proof}

\begin{lem}
\label{Brauer-Speiser} If $\chi{}^{}_{G, V}$ is real
valued, then ${\rm Schur}{}^{}_{G, V}\leqslant 2$.
\end{lem}

\begin{proof}  This is the result of R.~Brauer and
A.~Speiser, \cite{Brauer} (see also \cite{Fields},
\cite{Benard}, and \cite[Corollary~2.4,
p.\,277]{Feit}). \quad $\square$
\renewcommand{\qed}{}\end{proof}

\begin{lem}
\label{rational} Suppose that $\chi{}^{}_{G, V}$ is
rational and ${\,\rm Schur}{}^{}_{G, V}\neq 1$. Then
$D$ is 
a quaternion $\Q$-algebra, $n=2r$, and exactly one of
the following two possibilities holds:

\begin{itemize}
\item[(i)] the $G$-module $V$ is orthogonal and $D$ is
  indefinite,  i.e., $D\otimes_{\Q}\R$ is
$\R$-iso\-mor\-phic to $\M_2(\R)$;

\item[(ii)] the $G$-module $V$ is symplectic and $D$ is
  definite,  i.e., $D\otimes_{\Q}\R$ is
$\R$-isomorphic to $\H$.
\end{itemize}
\end{lem}

\begin{proof}
Since $\chi{}^{}_{G, V}$ is rational and ${\rm
Schur}{}^{}_{G, V}\neq 1$, Lemmas \ref{Dtrace} and
Lemma~\ref{Brauer-Speiser} imply that $Z_G=\Q$ and
${\rm Schur}{}^{}_{G, V}=2$. By \eqref{index}, this
implies that $\dim_{\Q}(D)=4$. Since $D$ is a division
$\Q$-algebra, the latter equality implies  that
 $D$ is a quaternion $\Q$-algebra (see,
 e.g.,\,\cite[\S 13.1]{PA}). Corollary \ref{nr} now
implies that $n^2=4r^2$, i.e.,
 $n=2r$.

 Since the $\R$-algebra $\Q G
 \otimes_{\Q}\R$ is
 simple, the natural
 surjection $\Q G
 \otimes_{\Q}\R\to\R G$ is
 an isomorphism of $\R$-algebras.
Clearly, the linear group $G$ is defined
 over $\R$ if and only if $\R G$
 is isomorphic to  $\M_n(\R)$. Taking
 into account the $\R$-algebra isomorphisms
 $$\Q G
 \otimes_{\Q}\R\simeq
 \M_{n/2}(D)\otimes_{\Q}\R\simeq\M_{n/2}(D\otimes_{\Q}\R),$$
 from this we deduce that $G$
 is defined
 over $\R$ if and only if
 $D$
  indefinite .
 Since a self-dual simple $G$-module
 is defined over $\R$ if and only if
 it is orthogonal
 (see, e.g., \cite[Section 13.2]{SR}),
 this completes the proof.
\quad $\square$
\renewcommand{\qed}{}\end{proof}

\begin{thm}\label{quater}
Suppose that $n$ is even and there exists a quaternion
$\Q$-algebra $H$ such that the $\Q$-algebras $\Q G$ and
$\M_{n/2}(H)$ are isomorphic. Then
\begin{itemize}
\item[(i)] for each imaginary quadratic subfield $F$ of
$H$, there exists a $G$-invariant lattice
$\Lambda^{(F)}$ in $V$ of rank $2n$ such that the
complex torus $V/\Lambda^{(F)}$ is isomorphic to a
product of mutually isogenous elliptic curves with
complex multiplication by an order  of $F$; \item[(ii)]
there exists a $G$-invariant lattice $\Lambda$ in $V$
of rank $2n$ such that the $\Q$-algebras
$\End(V/\Lambda)\otimes\Q$ and $\M_{n/2}(H)$ are
isomorphic.
\end{itemize}
\end{thm}

\begin{proof} Since the center of $\M_{n/2}(H)$ is $\Q$
(embedded by $a\mapsto a\cdot I_{n/2}$), we have
\begin{equation}\label{zq} Z_G=\Q.
\end{equation}

Fix an isomorphism of $\Q$-algebras
$$\tau\colon \Q G \longrightarrow \M_{n/2}(H).$$
The $\Q$-linear space $H^{n/2}$ carries a natural
structure of {\it left} $\M_{n/2}(H)$-module. Fix in
$H^{n/2}$ a nonzero $\tau(\Z G)$-stable finitely
generated additive subgroup $\Pi$ (for instance, take
$\Pi=\tau(\Z G)\cdot v$ for a nonzero vector $v\in
H^{n/2}$). Since $H^{n/2}$ is a $\Q$-linear space,
$\Pi$ is a free abelian group of finite rank and every
its basis consists of linearly independent elements of
$H^{n/2}$. Since the linear subspace $\Q \Pi$ of
$H^{n/2}$ is stable with respect to $\Q \tau(\Z G)=\tau
(\Q G)=\M_{n/2}(H)$, we have $\Q \Pi=H^{n/2}$. Hence
the natural map is an isomorphism of $\Q$-linear spaces
\begin{equation}\label{pi}
\Pi\otimes\Q\overset{\simeq}\longrightarrow H^{n/2}.
\end{equation}

The $\Q$-linear space  $ H^{n/2}$ carries also a
natural structure of {\it right} $H$-module. The
structures of left $\M_{n/2}(H)$-module and right
$H$-module on $ H^{n/2}$ yield the $\Q$-linear
embeddings of respectively $\M_{n/2}(H)$ and $H$ into
$\End_{\Q}(H^{n/2})$. By Wedderburn's theorem (see,
e.g.,\,\cite[Ch. XVII, \S 3, Corollary 3]{Lang}), the
images of these embeddings are the centralizers of
each other in $\End_{\Q}(H^{n/2})$.

Put $\Pi_{\R}:=\Pi\otimes\R$ and $H_{\R}:=
H\otimes_{\Q}\R$. We naturally identify $\Pi$ with
$\Pi\otimes 1$ in $\Pi_{\R}$, and $H$ with $H\otimes 1$
in $H_{\R}$. By \eqref{pi}, the following isomorphisms
of $\R$-linear spaces hold:
\begin{equation}\label{otimes}
\Pi_{\R} \simeq H^{n/2}\otimes_{\Q}\R\simeq
H_{\R}^{n/2}.
\end{equation}
From  \eqref{otimes} we deduce that $\Pi_{\R}$ is a
$2n$-dimensional $\R$-linear space that carries the
natural structures of left $\M_{n/2}(H_{\R})$-module
and right $H_{\R}$-module, and $\Pi$ is a {\it lattice}
 in $\Pi_{\R}$ of rank $2n$.
These structures yield the $\R$-linear embeddings
$$\iota_{l}\colon\M_{n/2}(H_{\R})
\hookrightarrow \End_{\R}(\Pi_{\R}),\quad
 \iota_{r}\colon H_{\R}\hookrightarrow
 \End_{\R}(\Pi_{\R})$$ whose images
 are
 the centralizers of each other in
 $\End_{\R}(\Pi_{\R})$.

Since the elements of $G$ are invertible in $\Q G$,
composition of the following embeddings of
$\Q$-algebras
\begin{equation}\label{chain}
\Q G \overset{\overset{\tau}\simeq} \longrightarrow
\M_{n/2}(H) \overset{\rm id}\hookrightarrow
\M_{n/2}(H_{\R})
\overset{\iota_{l}}\hookrightarrow\End_{\R}(\Pi_{\R})
\end{equation}
embeds $G$ into the group of invertible elements of
$\End_{\R}(\Pi_{\R})$. This defines an $\R$-linear
action of $G$ on $\Pi_{\R}$. By construction, the
lattice $\Pi$ in $\Pi_{\R}$ is $G$-invariant.

We want to define on the real linear space $\Pi_{\R}$ a
structure of complex linear space in such a way that
the algebra of its $\C$-linear transformations contains
$\iota_l(\M_{n/2}(H_{\R}))$. In order to do this,
choose an element $c \in H_{\R}$ with $c^2=-1$: using
that $\R$-algebra $H_\R$ is isomorphic to either $\H$
or $\M_2(\R)$ (see, e.g., \cite[\S\S 1.7, 13.2]{PA}),
it is easy to see that such $c$ exists and, moreover,
since $(aca^{-1})^2=-1$ for every invertible element
$a$ of $H_\R$, the set of such $c$
is uncountable (the latter fact will be used below).
Define now the complex structure on $\Pi_{\R}$ by
letting $\iota_r(c)$ be the multiplication by~$i$. Let
$V_c$ be the $n$-dimensional complex linear space
defined by this complex structure on $\Pi_{\R}$. Since
$\iota_r(c)$ commutes with the elements of
$\iota_l(\M_{n/2}(H_{\R}))$, we have
\begin{equation}\label{inclc}
\iota_{l}(\M_{n/2}(H_{\R})) \subset \End_{\C}(V_c).
\end{equation}
From \eqref{chain} and \eqref{inclc} we deduce that the
action of $G$ on $\Pi_\R$ defined above is $\C$-linear
with respect to this complex structure. Also,
\eqref{chain} and \eqref{inclc} clearly yield a nonzero
homomorphism of $\C$-algebras $\Q G\otimes_\Q\C \to
\End_{\C}(V_c)$ that endows $V_c$ with a nonzero
structure of 
$\Q G\otimes_{\Q}\C$-module of $\C$-dimension $n$. On
the other hand, $V$ is a nontrivial $\Q
G\otimes_{\Q}\C$-module of $\C$-dimension $n$ as well.
Notice now that there is  a unique (up to isomorphism)
$\Q G\otimes_{\Q}\C$-module of $\C$-dimension $n$
since, by \eqref{zq} and Lemma~\ref{i}, the
$\C$-algebra $\Q G\otimes_{\Q}\C$ is isomorphic to
$\M_n(\C)$. This
 implies that there is an isomorphism of
 $\Q G\otimes_{\Q}\C$-modules
\begin{equation}\label{isoac}
\nu_c\colon V_c\longrightarrow V.
\end{equation}
It follows from $G \subset \Q G\subset \Q
G\otimes_{\Q}\C$ that $\nu_c$ is an isomorphism of
$G$-modules.

Determine now the structure of the endomorphism algebra
of complex torus $V_c/\Pi$, i.e., the $\Q$-algebra
$$\End^0(V_c/\Pi):=\End(V_c/\Pi)\otimes \Q.$$
 Recall from \cite{OZ} that the {\it Hodge algebra}
$\HDG(V_c/\Pi)$ of $V_c/\Pi$ is the smallest
$\Q$-subalgebra $B$ of $\End_{\Q}(\Pi\otimes
{\Q})\subset \End_{\R}(\Pi_{\R})$ such that
$\iota_r(c)\in \R B $ and $\End^0(V_c/\Pi)$ coincides
with the centralizer of $B$ in $\End_{\Q}(\Pi\otimes
{\Q})$. Since $c\in H_{\R}$, we conclude that
$\HDG(V_c/\Pi)\subset \iota_r(H)$. Since the
centralizer of $\iota_r(H)$ in $\End_{\Q}(\Pi\otimes
{\Q})$ is $\iota_l(\M_{n/2}(H))$, this implies that
$\End^0(V_c/\Pi)\supset \iota_l(\M_{n/2}(H))$. As
$c\notin \R\cdot 1$, we have
$\dim_{\Q}\HDG(V_c/\Pi)\geqslant 2$. Hence, since $H$
is a quaternion $\Q$-algebra, $\HDG(V_c/\Pi)$ is either
$\iota_r(H)$ or $\iota_r(F)$ where $F$ is a quadratic
subfield of $H$. If $\HDG(V_c/\Pi)=\iota_r(F)$, then by
dimension reason,
\begin{equation}\label{fC} F\otimes_\Q \R=\R \cdot 1
+\R \cdot c \simeq \C, \end{equation}
 and
therefore $F$ is an imaginary quadratic field. If
$\HDG(V_c/\Pi)=\iota_r(H)$, then $\End^0(V_c/\Pi)$ is
the centralizator of $\iota_r(H)$, i.e.,
$\End^0(V_c/\Pi)=\iota_l(\M_{n/2}(H))$.

Prove now that when $c$ varies, all  possibilities for
$\HDG(V_c/\Pi)$ do occur, i.e.,
\begin{enumerate} \item[(a)] if $F$ is an imaginary
quadratic subfield of $H$, then for some $c$,
\begin{equation}\label{cF}
\iota_r(F)=\HDG(V_c/\Pi);
\end{equation}
 \item[(b)] there exists $c$ such that
$\End^0(V_{c}/\Pi)=\iota_l(\M_{n/2}(H))$.
\end{enumerate}

First, if $F$ is as in (a), then $F\otimes_{\Q}
\R\simeq \C$ and therefore there is an element $c_F\in
F\otimes_\Q{\R}\subset H_\R$ such that $c_F^2=-1$ (in
fact, there are exactly two such elements). Clearly,
then \eqref{cF} holds for $c=c_F$. This proves (a).

Second, notice that clearly the set of imaginary
quadratic subfields of $H$ is at most countable. For
every such field $F$, the intersection of
$F\otimes_\Q\R\simeq \C$ with the set $S:=\{c\in
H_{\R}\mid c^2=-1\}$ consists of two elements. Since
$S$ is uncountable,
this implies that there exists $c_0 \in H_{\R}$ such
that $c_0^2=-1$ and $c_0$ does not lie in
$F\otimes_{\Q}\R$ for every imaginary quadratic
subfield $F$ of $H$. Hence
$\HDG(V_{c_0}/\Pi)=\iota_r(H)$. This proves (b).

Determine now the structure of $V/\Pi_c$ in case when
$\HDG(V_c/\Pi)=\iota_r(F)$ where $F$ is an imaginary
quadratic subfield of $H$. The definition of Hodge
algebra implies that $\iota_r(F)\subset
\End^0(V_c/\Pi)$, i.e., $\Q\Pi$ is $\iota_r(F)$-stable.
From \eqref{fC} and the definition of $V_c$ we deduce
that
$$
\iota_r(F\otimes_\Q\R)=\C\cdot \I_{V_c}.$$ Since
$\Q\Pi$ is $\iota_r(F)$-stable, there is an order
$\mathcal O'$ in $F$ such that $\Pi$ is
$\iota_r(\mathcal O')$-stable. This endows $\Pi$ with a
structure of $\mathcal O'$-module. By a theorem of
Z.~I.~Borevich and D.~K.~Faddeev \cite{BF1} (see also
\cite{BF2}, \cite{BF3}, \cite[Satz 2.3]{Schoen}),
this
$\mathcal O'$-module splits into a direct sum of $n$
submodules of rank~$1$,
\begin{equation}\label{spl} \Pi=\Gamma_1\oplus
\ldots \oplus \Gamma_n.
\end{equation}
Clearly, each $\Gamma_j\otimes\R$ is a one-dimensional
$\C$-linear space and \eqref{spl} implies that
$$\textstyle
V_c/\Pi\quad \mbox{is isomorphic to}\quad\prod_{j=1}^n
(\Gamma_j\otimes\R)/\Gamma_j.$$ Every
$(\Gamma_j\otimes\R)/\Gamma_j$ is an elliptic curve
with complex multiplication by $\mathcal O'$.
Therefore these curves are mutually isogenous.

To complete the proof, it only remains to remark that
due to the existence of isomorphism \eqref{isoac}, the
complex tori $V_c/\Pi$ and $V/\nu_c(\Pi)$ are
isomorphic. So in the above cases (a) and (b), putting
respectively $\nu_c(\Pi):=\Lambda^{(F)}$ and $\Lambda$,
we obtain respectively the proofs of statements (i) and
(ii) of the theorem. \quad $\square$
\renewcommand{\qed}{}\end{proof}

We can now give a criterion of the existence of a
nonzero $G$-invariant lattice.

\begin{thm}\label{criterion}
{\rm (A)} The following properties are equivalent:
\begin{enumerate}
\item[\rm(a)] there is a nonzero $G$-invariant lattice
in $V$; \item[\rm(b)] there is a $G$-invariant lattice
in $V$ of rank $2n$; \item[\rm(c)] one of the following
conditions hold:
\begin{enumerate}
\item[\rm(i)] ${\rm Schur}^{}_{G,V}=1$ and
$\chi{}^{}_{G, V}$ is either rational or imaginary
quadratic; \item[\rm(ii)] ${\rm Schur}^{}_{G,V}=2$ and
$\chi{}^{}_{G, V}$ is rational.
\end{enumerate}
\end{enumerate}

{\rm (B)} A $G$-invariant lattice $\Lambda$ in $V$ of
rank $n$ exists if and only if $G$ is defined over
$\Q$, i.e., ${\rm Schur}{}_{G,V}=1$ and $\chi{}^{}_{G,
V}$ is rational. For such $\Lambda$ and every nonreal
$c\in \C$, the additive subgroup $\Lambda+c\Lambda$ of
$V$ is a $G$-invariant lattice in $V$ of rank $2n$.
\end{thm}
\begin{proof}
(A) The equivalence of (a) and (b) follows from (B)
that is proved below.

Assume now that (a) holds. If $\chi{}^{}_{G, V}$ is not
rational, then Lemmas \ref{Dtrace} and
\ref{nonrational} imply that ${\rm Schur}^{}_{G,V}=1$
and $\chi{}^{}_{G, V}$ is imaginary quadratic. If
$\chi{}^{}_{G, V}$ is rational, then Lemma
\ref{Brauer-Speiser} yields ${\rm
Schur}^{}_{G,V}\leqslant 2$. This proves that (a)
implies (c).

Conversely, assume that (c) holds. If (ii) is
fulfilled, then \eqref{D}, Lemma~\ref{rational}, and
Theorem~\ref{quater} imply that (a) holds. Consider now
the case when (i) is fulfilled.

If $\chi{}^{}_{G, V}$ is rational, then by definition
of the Schur index, $G$ is defined over $\Q$. It is
known (see, e.g.,\,\cite[(73.5)]{CR}) that then $G$ is
defined over $\Z$ , i.e., there exists a basis
$e_1,\ldots, e_n$ in $V$ such that the lattice
\begin{equation}\label{Zn}
\Z e_1+\ldots+\Z e_n \end{equation} is $G$-invariant.
Thus in this case (a) holds as well.

Finally, assume that $\chi{}^{}_{G, V}$ is imaginary
quadratic. Since ${\rm Schur}^{}_{G,V}=1$, there exists
a $G$-invariant $\Q (\chi{}^{}_{G, V})$-form $L$ of
$V$. Let $\mathcal O$ be the maximal order of $\Q
(\chi{}^{}_{G, V})$. Take any nonzero vector $v\in L$
and let $\Lambda$ be the submodule of $\mathcal
O$-module $L$ generated by the $G$-orbit of $v$,
\begin{equation}\label{O}
\textstyle \Lambda:=\sum_{g\in G}{\mathcal O} g(v).
\end{equation}
Since $\mathcal O$ is a Dedekind ring (see,
e.g.,\,\cite[\S 18]{CR}) and $\Lambda$ is a finitely
generated torsion free $\mathcal O$-module, the latter
is isomorphic to a direct sum of some fractional ideals
${\mathcal I}_1,\ldots, {\mathcal I}_d$ of $\Q
(\chi{}^{}_{G, V})$ (see, e.g.,\,\cite[(22.5)]{CR}).
Hence there are linearly independent over $\mathcal O$
vectors $v_1,\ldots, v_d\in L$ such that
\begin{equation}\label{directt}
\Lambda={\mathcal I}_1v_1+\ldots +{\mathcal I}_d v_d.
\end{equation}
Since the fraction field of $\mathcal O$ is $\Q
(\chi{}^{}_{G, V})$, vectors $v_1,\ldots, v_d$ are
linearly independent over $\Q (\chi{}^{}_{G, V})$ as
well, and since $L$ is a $\Q (\chi{}^{}_{G, V})$-form
of $V$, they are linearly independent over $\C$. Notice
now that since $\Q (\chi{}^{}_{G, V})$ is an imaginary
quadratic number field, all its fractional ideals are
lattices (of rank 2) in $\C$. This and \eqref{directt}
imply now that $\Lambda$ is a nonzero lattice in $V$.
On the other hand, \eqref{O} clearly implies that
$\Lambda$ is $G$-invariant. Hence (a) holds. This
completes the proof that (c) implies (a).

(B) We have already proven that if ${\rm
Schur}^{}_{G,V}=1$ and $\chi{}^{}_{G, V}$ is rational,
then there exists a $G$-invariant lattice of rank $n$,
namely, lattice~\eqref{Zn}. Conversely, let $\Lambda$
be a $G$-invariant lattice of rank $n$. Since
$\C\Lambda=V$ because of the irreducibility of $G$, the
equality $\rk(\Lambda)=\dim_{\C}V$ implies that every
basis $e_1,\ldots, e_n$ of the $\Z$-module $\Lambda$ is
a basis of the $\C$-linear space $V$. Hence $\Q
\Lambda$ is a $G$-invariant $\Q$-form of $V$, i.e., $G$
is defined over~$\Q$. Therefore in this case we have
\begin{equation}\label{ds}\Lambda+c\Lambda= (\Z
+c\Z)e_1\oplus\ldots \oplus(\Z  +c\Z) e_n.
\end{equation}
The condition $c\notin \R$ implies that $\Z +c\Z$ is a
lattice of rank $2$ in $\C$, and then from \eqref{ds}
we deduce that $\Lambda +c\Lambda$
 is a lattice of rank $2n$, which is clearly
$G$-invariant. This completes the proof.
 \quad
$\square$
\renewcommand{\qed}{}\end{proof}

\begin{cor}\label{refll} If ${\rm Schur}
{}^{}_{G, V}\geqslant 3$, then there are no nonzero
$G$-invariant lattices~in~$V$.
\end{cor}

\begin{cor}\label{nod} Suppose that the greatest common
divisor of the integers $\ker(u)$, where $u$ runs
through $\Q G$, is equal to $1$. Then
\begin{enumerate}
\item[\rm(i)] a nonzero $G$-invariant lattice in $V$
exists if and only if $\chi{}^{}_{G, V}$ is either
rational or imaginary quadratic; \item[\rm(ii)] a
$G$-invariant lattice in $V$ of rank $n$ exists if and
only if $\chi{}^{}_{G, V}$ is rational.
\end{enumerate}
\end{cor}

\begin{proof} This follows from Lemma \ref{split} and
Theorem \ref{criterion}. \quad $\square$
\renewcommand{\qed}{}\end{proof}

\begin{lem} \label{av}Let $\Lambda$ and $\Lambda'$ be
lattices of rank $2n$ in $V$ such that
$\Lambda'\subseteq \Lambda$.
\begin{enumerate} \item[\rm(a)]The
following properties are equivalent:
\begin{enumerate}
\item[\rm(i)] $V/\Lambda$ is an abelian variety;
\item[\rm(ii)] $V/\Lambda'$ is an abelian variety.
\end{enumerate}
\item[\rm(b)] If {\rm(i)} and {\rm(ii)} hold, then the
abelian varieties $V/\Lambda$ and $V/\Lambda'$ are
isogenous.
\end{enumerate}
\end{lem}
\begin{proof} (a) Let (i) holds. This means
 that $V/\Lambda$ admits a polarization $\Psi$,
 i.e., there exists a positive-definite Hermitian form
$\Psi\colon V\times V\to \C$ such that its imaginary
part assumes integer values on $\Lambda\times \Lambda$
(see, e.g.\,\cite{MumfordAV}). Since
$\Lambda'\subseteq\Lambda$, the same $\Psi$ is a
polarization for $V/\Lambda'$; whence~(ii).

Conversely, let (ii) holds and let $\Psi'$ be a
polarization for $V/\Lambda'$. Since ${\rm
rk}\,(\Lambda)={\rm rk}\,(\Lambda')$, we have
\begin{equation}\label{indexl}
[\Lambda:\Lambda']<\infty.
\end{equation}  Then clearly,
$[\Lambda:\Lambda']^2\cdot \Psi'$ is a polarization for
$V/\Lambda$; whence (i).

(b) The claim follows from inequality \eqref{indexl}
meaning that $V/\Lambda$ is the quotient of
$V/\Lambda'$ by a finite subgroup.
 \quad
$\square$
\renewcommand{\qed}{}\end{proof}

\begin{lem}\label{Katsura}
If an abelian variety is isogenous to a self-product of
an elliptic curve with complex multiplication, then, in
fact, it is isomorphic to a product of mutually
isogenous elliptic curves with complex multiplication.
\end{lem}
\begin{proof} This is proved in \cite{Katsura} (see
also \cite{Shioda}, \cite{Schoen}). \quad $\square$
\renewcommand{\qed}{}\end{proof}
\section{\bf
The case
of \boldmath${\rm Schur}{}^{}_{G, V}=1$} The following
theorem gives a description of $V/\Lambda$ in case when
${\rm Schur}{}^{}_{G, V}=1$.

\begin{thm}\label{main2}
Suppose that ${\rm Schur}{}^{}_{G, V}=1$. If there
exists a nonzero $G$-invariant lattice $\Lambda$ in
$V$, then the following properties hold:
\begin{enumerate}
\item[(i)] $Z_{G}$ is either $\Q$ or an imaginary
quadratic number field. \item[(ii)] ${\rm
rk}\,(\Lambda)=n$ or $2n$. \item[(iii)] If ${\rm
rk}\,(\Lambda)=n$, then $Z_{G}=\Q$ and $G$ is defined
over $\Q$.
 \item[(iv)] Suppose that
 ${\rm rk}\,(\Lambda)=2n$.
Let $\O$ be the maximal order in $Z_{G}$. Fix a
$Z_{G}$-algebra isomorphism $($existing by
{\rm\eqref{<=>}}$)$
 $$\tau:\M_n(Z_{G})\overset{\simeq}\longrightarrow \Q G.$$
 Then there is a
 lattice $\Lambda'$
in $V$ that enjoys the following
properties:
\begin{list}{}{\itemindent=0mm \leftmargin=10mm}
\item[\rm (iv${}_1$)] $\Lambda'\supseteq \Lambda$;
\item[\rm (iv${}_2$)] $\Lambda'$ is
$\tau(\M_n(\O))$-invariant; \item[\rm (iv${}_3$)] there
exists a lattice $\Gamma$ of rank $2$ in $\C$ and a
$\C$-linear isomorphism $\nu:\C^n
\overset{\simeq}\longrightarrow V$ such that
$\O\Gamma=\Gamma$ and $\nu(\Gamma^n)=\Lambda'$;
\item[\rm (iv${}_4$)] the complex torus $V/\Lambda$ is
an abelian variety isogenous to a self-product of an
elliptic curve; \item[\rm (iv${}_5$)] if $Z_{G}$ is an
imaginary quadratic number field, then $V/\Lambda$ is
isomorphic to a product of mutually isogenous elliptic
curves with complex multiplication  by $\mathcal O$.
\end{list}
\end{enumerate}
\end{thm}

\begin{proof} Lemmas~\ref{nonrational0} and \ref{Dtrace}
imply (i), Lemma~\ref{n2n} implies (ii), and
Theorem~\ref{criterion} and Lemma~\ref{Dtrace} imply
(iii).  Assume now that ${\rm rk}(\Lambda)=2n$ and
prove (iv).

 Clearly, the map $\O\otimes\Q \to Z_{G},\
a\otimes r\mapsto a r$, is a ring isomorphism that
yields the ring isomorphism
$$\M_n(\O)\otimes\Q \overset{\simeq}
\longrightarrow \M_n(Z_{G}),\ z\otimes r\mapsto z r.$$
Thus $\M_n(\O)$ is an order in the $\Q$-algebra
$\M_n(Z_G)$; whence $\tau(\M_n(\O))$ is an order in the
$\Q$-algebra $\Q G$. Since $\Z G$ is an order in $\Q G$
as well, $\Z G\cap \tau(\M_n(\O))$ is a subgroup of
finite index in $\tau(\M_n(\O))$. Since $\Lambda$ is
$\Z G$-invariant, this entails that there are only
finitely many sets of the form $z(\Lambda)$, where
$z\in \tau(\M_n(\O))$. Every such set $z(\Lambda)$ is a
finitely generated additive subgroup in $V$, and
$[z(\Lambda):(\Lambda\cap z(\Lambda))]<\infty$. This
implies that the sum of these subgroups,
\begin{equation}\label{lambda'}
\textstyle \Lambda':=\sum_{z\in\tau(\M_n(\O))}
z(\Lambda), \end{equation} is a
$\tau(\M_n(\O))$-invariant lattice in $V$ containing
$\Lambda$ as a subgroup of finite index.

The $\tau(\M_n(\O))$-module $\Lambda'$ is faithful.
Indeed, notice that since $\Lambda'$ is
$\tau(\M_n(\O))$-invariant and $\tau(\M_n(\O))$ is an
order in $\Q G$, the $\Q$-linear subspace $\Q\Lambda'$
in $V$ is $\Q G$-invariant. Therefore if $z\Lambda'=0$
for $z\in \tau(\M_n(\O))$, then $z$ lies in the kernel
of natural $\Q$-algebra homomorphism $\Q G \to
\End_{\Q}(\Q\Lambda')$. Since $\Q G$ is a simple
$\Q$-algebra, this kernel is trivial. Thus $z=0$;
whence the faithfulness.

By construction, (iv${}_1$) and (iv${}_2$) hold. We are
going to prove that (iv${}_3$), (iv${}_4$), and
(iv${}_5$) hold as well.

Tensoring $\tau$ by $\C$ over $Z_{G}$, we obtain a
$\C$-algebra isomorphism
\begin{equation}\label{iso1}
 \M_n(\C)=\M_n(Z_{G})\otimes_{Z_{G}}\!\C
\overset{\simeq}\longrightarrow \Q
G\otimes_{Z_{G}}\!\C.
\end{equation}

On the other hand, by Lemma~\ref{i}, we have the
$\C$-algebra isomorphism \eqref{psi} Composing
isomorphisms \eqref{iso1} and \eqref{psi}, we get a
$\C$-algebra isomorphism
$$\tau_{\C}\colon \M_n(\C)
\overset{\simeq}\longrightarrow
 \End_{\C}(V).$$

Consider the coordinate $\C$-linear space $\C^n$
endowed with the natural structure of left
$\M_n(\C)$-module. Since $\C^n$ is the unique (up to
isomorphism) left $\M_n(\C)$-module of $\C$-dimension
$n$, there is a $\C$-linear isomorphism $\varphi\colon
\C^n \overset{\simeq}\longrightarrow V$ such that
$$\varphi(y(v))=\tau_{\C}(y)
(\varphi(v)) \hskip 3mm \text{for all} \ v\in \C^n, y
\in \M_n(\C).$$ Clearly,
$$\Lambda'':=\varphi^{-1}(\Lambda')$$
is an $\M_n(\O)$-invariant lattice of rank $2n$ in
$\C^n$. Let $e_1,\ldots , e_n$ be the standard basis in
$\C^n$. Since $\Lambda''$ is $\M_n(\O)$-invariant, it
is easily seen that there exists a lattice $\Gamma$ in
$\C$ such that
\begin{equation}\label{gamman}
\O\Gamma=\Gamma\hskip 2mm \text{and}\hskip 2mm
\Lambda''=\Gamma e_1+\ldots+\Gamma e_n.
\end{equation} Since ${\rm rk}(\Lambda'')=2n$, we
deduce from the second equality in \eqref{gamman} that
${\rm rk}(\Gamma)=2$, and the complex torus
$V/\Lambda'$ is isomorphic to the self-product of
elliptic curve $\C/\Gamma$, hence is an abelian
variety. Since $\Lambda$ is a subgroup of finite index
in $\Lambda'$, Lemma~\ref{av} implies that $V/\Lambda$
is an abelian variety as well and $V/\Lambda$ and
$V/\Lambda'$ are isogenous. This proves~(iv${}_3$) and
(iv${}_4$).

Now assume that $Z_{G}$ is an imaginary quadratic
field. Then $\mathcal O\neq \Z$, and the first equality
in \eqref{gamman} yields that the elliptic curve
$\C/\Gamma$ has complex multiplication. Hence
$V/\Lambda$ is isogenous to a self-product of an
elliptic curve with complex multiplication. Now
(iv${}_5$) follows from Lemma~\ref{Katsura}. \quad
$\square$
\renewcommand{\qed}{}\end{proof}

Combining Theorem~\ref{main2} with the results of
Section 2, we obtain the following applications.

\begin{thm}
\label{nonrationalT} Suppose that $\chi{}^{}_{G, V}$ is
not rational. If there exists a nonzero
$G$-in\-va\-ri\-ant lattice $\Lambda$ in $V$, then
\begin{itemize}
\item[\rm (i)] $Z_{G}$ is an imagi\-nary quadratic
number field; \item[\rm (ii)] $\rk(\Lambda)=2n$;
\item[\rm (iii)] $V/\Lambda$ is isomorphic to a product
of mutually isogenous elliptic curves with complex
multiplication by
$Z_G$.
\end{itemize}
\end{thm}

\begin{proof}
This follows from the combination of Lemma
\ref{nonrational} and Theorem \ref{main2}. \quad
$\square$
\renewcommand{\qed}{}\end{proof}

\begin{cor}
[E.~B.~Vinberg]
\label{vbgT} Suppose that the
$G$-module $V$ is not self-dual. If there exists a
nonzero $G$-invariant lattice $\Lambda$ in $V$, then
$\rk(\Lambda)=2n$ and $V/\Lambda$ is isogenous to a
self-product of an elliptic curve with complex
multiplication.
\end{cor}

\begin{proof}
This immediately follows  from Theorem
\ref{nonrationalT}. \quad $\square$
\renewcommand{\qed}{}\end{proof}

\begin{thm} \label{gcd}
Suppose that the greatest common divisor of the
integers $\ker(u)$, where $u$ runs through $\Q G$, is
equal to $1$. If there exists a nonzero $G$-invariant
lat\-ti\-ce $\Lambda$ in $V$, then the conclusions of
Theorem~{\rm \ref{main2}} hold true.
\end{thm}
\begin{proof}
This follows from the combination of Lemma~\ref{split}
and Theorem~{\rm \ref{main2}}. $\square$
\renewcommand{\qed}{}\end{proof}

\begin{rem}\label{refl}
Since $\dim_{\C}\ker(0)=n$ and
$\dim_{\C}(\ker(\I_V-r))=n-1$ for any reflection $r\in
\GL(V)$,  the condition that the greatest common
divisor of the integers $\ker(u)$, where $u$ runs
through $\Q G$, is equal to $1$, always holds for every
reflection group $G$. By Theorem \ref{gcd}, this
implies that the conclusion of Theorem \ref{main2}
holds true for every reflection group $G$ that admits a
nonzero invariant lattice.
\end{rem}

\begin{ex}\label{double}
Suppose that $n$ is a positive integer, $S$ is a $(n+1)$-element
set, $G$ is a doubly transitive permutation group of $S$ and $V$
is the $n$-dimensional $\C$-vector space of complex-valued
functions $f:S \to \C$ with $\sum_{s\in S}f(s)=0$. Clearly, the
natural linear representation of $G$ in $V$ is faithful and
defined over $\Q$; in particular, its character is rational. It is
well-known that the $G$-module $V$ is (absolutely) simple (see,
for instance, \cite[Sect. 2.3, Ex. 2]{SR}). Clearly, the Schur
index is $1$. It follows from Theorem \ref{criterion} that there
exist $G$-invariant lattices in $V$ of rank $n$ and $2n$. It
follows from Theorem \ref{main2} that if $\Lambda$ is a
$G$-invariant lattice of rank $2n$ in $V$ then $V/\Lambda$ is
isogenous to a self-product of an elliptic curve.
\end{ex}

\begin{ex} \label{L35}
Let $G$ be the simple group $L_3(5):=\PSL_3(\F_5)$. Then there
exists a simple complex $G$-module $V$ such that
$\dim_{\C}(V)=124$, $Z_G=\Q(\chi^{}_{G, V})=\Q(\sqrt{-1})$ and
 ${\rm Schur}^{}_{G, V}=1$,  see \cite[p.\,283]{Feit}. By
Theorem~\ref{criterion}, there are $G$-invariant lattices of rank
$248$ in $V$. It follows from Theorem \ref{main2} that if
$\Lambda$ is a $G$-invariant lattice of rank $248$ in $V$ then
$V/\Lambda$ is  isomorphic to a product of mutually isogenous
elliptic curves with complex multiplication by $\Q(\sqrt{-1})$.
\end{ex}

\begin{ex}
Let $p$ be an odd prime that is congruent to $3$ modulo $4$, $r$ a
positive integer and $q=p^{2r-1}$. Let $G$ be the  group
$\SL_2(\F_q)$. Then there exists a faithful
 simple complex $G$-module $V$ such that $\dim_{\C}(V)=(q-1)/2$,
$Z_G=\Q(\chi^{}_{G, V})=\Q(\sqrt{-q})=\Q(\sqrt{-p})$ and
 ${\rm Schur}^{}_{G, V}=1$,  see \cite[p.\, 4]{J}, \cite[p.\,284--285]{Feit}.
 By
Theorem~\ref{criterion}, there are $G$-invariant lattices of rank
$q-1$ in $V$. It follows from Theorem \ref{main2} that if
$\Lambda$ is a $G$-invariant lattice of rank $q-1$ in $V$ then
$V/\Lambda$ is isomorphic to a product of mutually isogenous
elliptic curves with complex multiplication by $\Q(\sqrt{-p})$.

\end{ex}

\section{\bf
The case
of  \boldmath${\rm Schur}{}^{}_{G, V}\neq 1$ }

The following theorem gives a description of
$V/\Lambda$ in case when ${\rm Schur}{}^{}_{G, V}\neq
1$ (i.e., in fact, when ${\rm Schur}{}^{}_{G, V}=2$, by
Theorem~\ref{criterion}).

\begin{thm}
\label{ratL} Suppose that ${\rm Schur}{}^{}_{G, V}\neq
1$. If there exists a nonzero $G$-invariant lattice
$\Lambda$ in $V$, then the following properties hold:
\begin{enumerate}
 \item[\rm (i)] $\chi^{}_{G, V}$ is
rational $($hence the $G$-module $V$ is
either orthogonal or symplectic$)$.
\item[\rm (ii)] $n$
is even.
 \item[\rm (iii)] There
 exists
  a quaternion $\Q$-algebra $H$ such that
the $\Q$-algebras $M_{n/2}(H)$ and $\Q G$ are
isomorphic.
\item[\rm (iv)] $\rk(\Lambda)=2n$.
\item[\rm (v)] 
Fix an order $\O$ in $H$. Then there exists a
two-dimensional complex torus $T$ that enjoys the
following properties: 
\begin{list}{}{\itemindent=0mm \leftmargin=10mm}
\item[${\rm (v_1)}$] $V/\Lambda$ is isogenous to a
self-product of $T$; \item[${\rm (v_2)}$] there exists
a ring embedding $\O \hookrightarrow \End(T)$;
\item[${\rm (v_3)}$] if the $G$-module $V$ is
orthogonal, then $H$ is   indefinite,  and $T$ and
$V/\Lambda$ are abelian varieties; \item[${\rm (v_4)}$]
if the $G$-module $V$ is symplectic, then $H$ is
definite, and $T$ and $V/\Lambda$ either are not
abelian varieties or
are isomorphic to the products of mutually isogenous
elliptic curves with complex multiplication.
\end{list}
\end{enumerate}
\end{thm}

\begin{proof}
Theorem~\ref{criterion}, Lemma~\ref{rational}, and
Lemma~\ref{n2n} imply (i), (ii), (iii), and (iv). Fix a
$\Q$-algebra isomorphism
$$\tau: \M_{n/2}(H)
\overset{\simeq}\longrightarrow \Q G$$ and an order
$\O$ in $H$. Both $\Z G$ and $\tau(\M_{n/2}(\O))$ are
the orders in $\Q$-algebra $\Q G$. The same argument as
in the part of proof of Theorem~\ref{main2} related to
formula \eqref{lambda'} shows that there are only
finitely many sets of the form $z(\Lambda)$, where
$z\in \tau(\M_{n/2}(\O))$, and the sum
$$\textstyle
\Lambda':= \sum_{z\in\tau(\M_{n/2}(\O))}
z(\Lambda)\subset V$$ is a
$\tau(\M_{n/2}(\O))$-invariant lattice in $V$
containing $\Lambda$ as a subgroup of finite index and
faithful as $\tau(\M_{n/2}(\O))$-module.

Tensoring $\tau$ by $\C$ over $\Q$ and composing with
\eqref{psi}, we get a $\C$-algebra isomorphism
\begin{equation*}
\tau_{\C}\colon \M_{n/2}(H)\otimes_{\Q}\C
\overset{\simeq}\longrightarrow \End_{\C}(V).
\end{equation*}
Since the $\C$-algebras $H\otimes_{\Q}\C$ and
$\M_2(\C)$ are isomorphic, we may (and shall) fix a
$\C$-algebra isomorphism
$$\kappa\colon
H\otimes_{\Q}\C \overset{\simeq}\longrightarrow
\M_2(\C).$$ It induces the $\C$-algebra isomorphism
$$\kappa_{n/2}: \M_{n/2}(H)
\otimes_{\Q}\C= \M_{n/2}(H\otimes_{\Q}\C)
\overset{\simeq}\longrightarrow \M_{n/2}(\M_2(\C)).
$$
Hence we obtain a $\C$-algebra isomorphism
$$\tau_{\C}\circ\kappa_{n/2}^{-1}
\colon \M_{n/2}(\M_2(\C))
\overset{\simeq}\longrightarrow \End_{\C}(V).$$

 Consider the coordinate $\C$-linear space $\C^n$
presented as the direct sum of $n/2$ copies of $\C^2$,
$$\C^n=(\C^2)^{n/2}=
\C^2\oplus \ldots \oplus \C^2,$$ and endowed with the
natural structure of left $\M_{n/2}(\M_2(\C))$-module.
Since $(\C^2)^{n/2}$ is the unique (up to isomorphism)
left $\M_{n/2}(\M_2(\C))(\simeq\!\M_n(\C))$-module of
$\C$-di\-men\-sion $n$, there is a $\C$-linear
isomorphism $ \varphi\colon
(\C^2)^{n/2}\overset{\simeq}\longrightarrow V$ such
that
$$
\varphi(y(v))=(\tau_{\C}\circ \kappa_{n/2}^{-1}) (y)
(\varphi(v)) \hskip 2mm \mbox{for all}\ v\in
(\C^2)^{n/2},\ y\in \M_{n/2}(\M_2(\C)).$$ We identify
$H$ with $H\otimes_{\Q} 1$. Then clearly,
$$\Lambda^{''}:=
\varphi^{-1}(\Lambda')$$ is a $\M_{n/2}(\kappa
(\mathcal O))$-invariant lattice of rank $2n$ in
$(\C^2)^{n/2}$. Using the $\M_{n/2}(\kappa (\mathcal
O))$-invariance of $\Lambda^{''}$, it is easy to see
that there exists a lattice $\Gamma$ in $\C^2$ such
that
\begin{equation}\label{deco}
\kappa(\O)\Gamma=\Gamma\hskip 2mm \mbox{and}\hskip 2mm
\Lambda^{''}=\Gamma^{n/2}= \Gamma \oplus \ldots \oplus
\Gamma.
\end{equation}
Since $\rk(\Lambda^{''})=\rk(\Lambda^{'})
=\rk(\Lambda)=2n$, we deduce from \eqref{deco} that
$\rk(\Gamma)= 4$. Hence $T:=\C^2/\Gamma$ is a
$2$-di\-men\-sional complex torus and we have the ring
embedding
\begin{equation}\label{kO}
\O\overset{\overset{\kappa}\simeq} \longrightarrow
\kappa(\O) \hookrightarrow \End(T).
\end{equation} Clearly, the complex
torus $(\C^2)^{n/2}/\Lambda^{''}
\simeq(\C^2/\Gamma)^{n/2}=T^{n/2}$ is isomorphic to
$V/\Lambda'$. Since $\Lambda$ is a subgroup of finite
index in $\Lambda'$, the complex torus $V/\Lambda$ is
isogenous to $T^{n/2}$. This proves ${\rm (v_1)}$
and~${\rm (v_2)}$.

Suppose that  $H$ is   indefinite,  i.e., by Lemma
\ref{rational}, the $G$-module $V$ is orthogonal. It
then follows from \cite[Theorem 4.3, p.\,152]{L} that
$T$ is an abelian surface. Now
 Lemma \ref{av} implies that
 $V/\Lambda$ is
 an abelian variety. This proves
~${\rm (v_3)}$.

Assume now that  $H$ is   indefinite,  i.e., by Lemma
\ref{rational}, the $G$-module $V$ is symplectic. By
\eqref{kO}, the $\Q$-algebra
$\End^0(T):=\End(T)\otimes\Q$ contains a subalgebra
isomorphic to $\O\otimes\Q\simeq H$ and therefore is
noncommutative.

Suppose that $T$ is an abelian surface. Using tables in
\cite{O}, one may then easily verify that $T$ is not
simple. On the other hand, if $T$ is isogenous to a
product of two non-isogenous elliptic curves, then
$\End^0(T)$ is commutative and cannot contain a
subalgebra isomorphic to $H$. If $T$ is isogenous to a
square of an elliptic curve without complex
multiplication, then $\End^0(T)$ is isomorphic to
$\M_2(\Q)$ and hence cannot contain such a subalgebra
as well. It follows that $T$ is isogenous to a square
of an elliptic curve with complex multiplication. This
implies that $V/\Lambda$ is isogenous to a self-product
of an elliptic curve with complex multiplication, and
hence, by Lemma~\ref{Katsura}, that $V/\Lambda$ is
isomorphic to a product of mutually isogenous elliptic
curves with complex multiplication.

Now assume that $T$ is not an abelian surface. Then
$T^{n/2}$ is not an abelian variety because every
polarization
on $T^{n/2}$ obviously induces a polarization
on (say, the first factor) $T$. It then follows from
Lemma \ref{av} that $V/\Lambda$ is not an abelian
variety as well. This proves ${\rm (v_4)}$. \quad
$\square$
\renewcommand{\qed}{}\end{proof}

The following examples show that both possibilities in
conclusion~${\rm (v_4)}$ of Theorem~\ref{ratL} may
indeed occur. In particular, there are finite
irreducible groups $G$ and $G$-invariant lattices
$\Lambda$ in $V$ such that ${\rm rk}(\Lambda)=2n$ and
the complex torus $V/\Lambda$ is {\it not} an abelian
variety.

\begin{ex} First, we can use that $H$ contains an
imaginary quadratic subfield $F$ (see below the proof
of Theorem~\ref{deform}). Then by Theorem~\ref{quater},
for $\Lambda=\Lambda^{(F)}$, the torus $V/\Lambda$ is
isomorphic to a product of mutually isogenous elliptic
curves with complex multiplication.

One can also give a more concrete example not
exploiting Theorem~\ref{quater}. Let $G$ be the image
of a (unique up to isomorphism) irreducible
$2$-dimensional complex representation of the
quaternion group. Fixing a basis in the representation
space $V$, we may (and shall) identify $V$ with $\C^2$
and $G$ with the matrix group
\begin{equation}\label{quatern}
{\fontsize{9pt}{3mm}\selectfont \Bigl\{
\pm\begin{bmatrix}
1&0\\
0&1
\end{bmatrix},
\pm\begin{bmatrix}
i&0\\
0&-i
\end{bmatrix},
\pm\begin{bmatrix}
0&1\\
-1& 0
\end{bmatrix},
\pm\begin{bmatrix}
0&i\\
i&0
\end{bmatrix}
\Bigr\}. }
\end{equation}
It is known, \cite[\S\,70]{CR}, that ${\rm
Schur}^{}_{G, V}=2$. The $G$-module $\C^2$ is
symplectic and \eqref{quatern} clearly yields that the
lattice
$\Lambda:=\{\bigl(\begin{smallmatrix}a\\
b\end{smallmatrix}\bigr)\in \C^2\mid a, b\in \Z+i\Z]\}$
is $G$-stable and $\C^2/\Lambda$ is isomorphic to the
square of the elliptic curve $\C/(\Z+i\Z)$ with complex
multiplication by $\Z[i]$.\quad $\square$
\end{ex}

\begin{ex}\label{non}
Consider the order $\Lambda:=\Z\1+\Z\i + \Z\j + \Z\k$
in the
Hamiltonian quaternion $\Q$-algebra
$H:=\bigl(\frac{-1, -1}{\Q}\bigr)=\Q\1+\Q\i + \Q\j +
\Q\k$. It is a
  lattice of rank $4$ in the underlying
  $4$-dimensional real
linear space $V$ of the quaternion $\R$-algebra $\H=
\R\1+\R\i + \R\j + \R\k$. The quaternion group
$G:=\{\pm\1, \pm\i, \pm\j, \pm\k\}$ acts $\R$-linearly
and faithfully on $V$ by {\it left} multiplication in
$\H$. This action is irreducible and $\Lambda$ is
$G$-invariant. Pick any ele\-ment $c\in \H$ such that
$c^2=-\1$, and endow $V$ with a structure of
$2$-dimensional complex linear space $V_c$ defining
multiplication by $i$ as {\it right} multiplication by
$c$ in~$\H$. Since left and right multiplications
commute, the action of $G$ on $V_c$ is $\C$-linear (and
irreducible). Thus we may (and shall) view $G$ as an
irreducible group of complex linear transformations of
$V_c$.

Assume now that $c\notin\R F$ for any imaginary
quadratic
 subfield $F$ of $H$
(such $c$'s do exist, see the proof of
Theorem~\ref{quater}). Consider the endomorphism ring
of $2$-dimensional complex torus $V/\Lambda$,
\begin{equation*}
\End(V/\Lambda):=\{u \in \End_{\C}(V)\mid
u(\Lambda)\subseteq \Lambda\}.
\end{equation*}
Our assumption on $c$ implies that every $\Q$-linear
endomorphism of $H$ that becomes an element of
$\End_{\C}(V)$
 after the extension of scalars
 from $\Q$ to $\R$ must
commute with right multiplication by every element of
$H$, and therefore is left multiplication by an
element of~$H$. Hence $\End(V/\Lambda)$ consists of
all $u \in H$ with $u\cdot\Lambda\subseteq\Lambda$.
It follows that $\End(V/\Lambda)$ coincides with the
set of left multiplications by elements of $\Lambda$
and therefore $\End(V/\Lambda)\simeq \Lambda$. Notice
that $H$ is a {\it definite} quaternion $\Q$-algebra.
Thus the endomorphism ring of two-dimensional complex
torus $V/\Lambda$ is an order in a definite
quaternion $\Q$-algebra. But it is known, \cite{O},
that there are no complex abelian surfaces whose
endomorphism ring is an order in a definite
quaternion $\Q$-algebra. This implies that
$V/\Lambda$ is {\it not} an abelian variety.

 Note that if $c\in\R F$ for an imaginary
quadratic subfield $F$ of $H$, then the endomorphism algebra of
$V/\Lambda$ is isomorphic to $H \otimes_{\Q}F \simeq {\rm M}_2(F)$
and therefore $V/\Lambda$ is isogenous to a square of an elliptic
curve with complex multiplication by an order of $F$. \quad
$\square$
\end{ex}

\begin{ex} Here is another example of the outcome (v${}_4$) of
Theorem~\ref{ratL}. Let $p$ be an odd prime, $r$ a positive
integer and $q=p^{2r}$.

Let $G$ be the  group $\SL_2(\F_q)$. Then there exists a faithful
 simple complex $G$-module $V$ such that $\dim_{\C}(V)=(q-1)/2$,
$\chi^{}_{G, V}$ is rational, ${\rm Schur}^{}_{G, V}=2$, and the
quaternion $\Q$-algebra $H$ from Theorem~\ref{ratL} (iii) is
ramified exactly at $p$ and $\infty$,  see \cite[p.\, 4]{J},
\cite[p.\,284--285]{Feit}. In particular, $H$ is definite. By
Theorem~\ref{criterion}, there are $G$-invariant lattices
$\Lambda$ of rank $q-1$ in $V$. \quad $\square$
\end{ex}

\begin{ex} Here is an example of the outcome (v${}_3$) of
Theorem~\ref{ratL}. Let $G$ be the simple group ${\rm HJ}$. Then
there exists a simple complex $G$-module $V$ such that
$\dim_{\C}(V)=336$, $\chi^{}_{G, V}$ is rational, ${\rm
Schur}^{}_{G, V}=2$, and the quaternion $\Q$-algebra $H$ from
Theorem~\ref{ratL} (iii) is
  indefinite,  see   \cite[p.\,283]{Feit}. By
Theorem~\ref{criterion}, there is a  $G$-invariant lattice
$\Lambda$ of rank $672$ in $V$. \quad $\square$
\end{ex}

According to Example~\ref{non}, in general there exist
$G$-invariant lattices $\Lambda$ such that the complex
torus $V/\Lambda$ {\it is not} an abelian variety.
However it appears that one can always replace
$\Lambda$ by another $G$-invariant lattice $\Delta$
such that $V/\Delta$ {\it is} an abelian variety. More
precisely, the following statement holds true.

\begin{thm} \label{deform} The following properties are equivalent.
\begin{enumerate}
\item[(i)] there exists a nonzero $G$-invariant lattice
in $V$; \item[(ii)] there exists a $G$-invariant
lattice $\Delta$ in $V$ such that $V/\Delta$ is
isomorphic to a product of mutually isogenous elliptic
curves with complex multiplication.
\end{enumerate}
\end{thm}
\begin{proof} Assume that (i) holds. Then
Theorems~\ref{criterion} and \ref{main2} reduce proving
(ii) to the cases when $\chi^{}_{G, V}$ is rational and
${\rm Schur}^{}_{G, V}$ is $1$ or $2$.

Consider the case when $\chi^{}_{G, V}$ is rational and ${\rm
Schur}^{}_{G, V}=1$. Let $\mathcal O$ be an order in a imaginary
quadratic number field. We have $\mathcal O=\Z+c\Z$ for some
non-real $c\in \C$. Theorem~\ref{criterion} then implies that
\eqref{ds} is a $G$-invariant lattice of rank $2n$; denote it by
$\Delta$. By construction, $\C/\mathcal O$ is an elliptic curve
with complex multiplication by $\mathcal O$, and \eqref{ds}
implies that $V/\Delta$ is isomorphic to $(\C/\mathcal O)^n$. Thus
in this case (ii) holds.

Consider now the case when $\chi^{}_{G, V}$ is rational and ${\rm
Schur}^{}_{G, V}=2$. Lemma~\ref{rational} and Theorem~\ref{quater}
then reduce proving (ii) to showing that every quaternion
$\Q$-algebra $H=\bigl(\frac{a, b}{\Q}\bigr)$ contains an imaginary
quadratic subfield. But the latter property indeed holds, since
the maximal subfields of $H$ are precisely (up to isomorphism) the
fields $\Q(\sqrt{ar_1^2+b r_2^2-ab r_3^2})$, where $r_1, r_2,
r_3\in \Q$ and $r_1^2+r_2^2+r_3^2\neq 0$ (see, e.g.,\,\cite[\S
13.1, Exercise 4]{PA}). This completes the proof. \quad $\square$
\renewcommand{\qed}{}\end{proof}

\section{\bf Quotients modulo
invariant lattices of reflection groups}

In this section we give another, geometric proof of
Theorem~\ref{main2} for reflection groups (regarding
the first proof see Remark~\ref{refl}). It provides a
more precise information on invariant lattices.

Recall that an element $r\in \GL(V)$ is called a ({\it
complex\,}) {\it reflection} if
\begin{enumerate}
\item[(i)] the order of $r$ is finite; \item[(ii)]
$\dim_{\C}({\rm ker}(\I_V-r))=n-1$.
\end{enumerate}
For such $r$, the linear subspace
\begin{equation}\label{l}
l_r:=(\I_V-r)(V)
\end{equation}
is one-dimensional, $r$-invariant, and $r$ acts on it
as scalar multiplication by a root of unity
$\theta_r\neq 1$.

\begin{rem}\label{rgr} This implies that the assumptions
and conclusions of Lem\-ma~\ref{split} and
Corollary~\ref{nod} hold if $G$ is a reflection group.
\end{rem}

 It is well known, \cite{ST} (see also
\cite{Cohen}, \cite{Popov}), that every finite
irreducible reflection group in $V$ is generated by $n$
or $n+1$ reflections and, in the last case, it contains
an irreducible reflection subgroup generated by $n$
reflections. Therefore describing invariant lattices of
finite irreducible reflection groups in $V$, we may
consider only the groups generated by $n$ reflections.
Let $G$ be such a group, and let $r_1,\ldots, r_n$ be a
system of reflections generating $G$. We put
$l_{r_j}=l_j$, $\theta_{r_j}=\theta_j$.

Since $G$ is finite, we may (and shall) fix a
$G$-invariant positive definite Hermitian inner product
$\langle\ {,}\ \rangle$ on $V$. For every reflection
$r\in G$, fix a vector $e_r\in l_r$ of length~$1$. We
then have
\begin{equation}\label{formula}
r(v)=v-(1-\theta_r)\langle v {,} e_r \rangle e_r,\quad
v\in V.
\end{equation}
We put $e_j:=e_{l_j}$.

Denote by $\mathcal L$ the set of all the lines $l_r$
where $r$
 runs through all the reflections in $G$. Since $G$ is
 irreducible, the mutual fixed point set of
 $r_1,\ldots,
 r_n$ is $\{0\}$; whence
\begin{equation}\label{direct}
V=l_1\oplus\ldots \oplus l_n.
\end{equation}

Let $\Lambda$ be a $G$-invariant lattice in $V$. We put
\begin{equation}\label{G}
\textstyle \Lambda^0:=\sum_{l\in \mathcal
L}\Lambda_l,\quad \text{where}\ \Lambda_l:=\Lambda\cap
l,
\end{equation}
and $\Lambda_j:=\Lambda_{l_j}$. Then $\Lambda^0$ is a
subgroup of $\Lambda$, hence $\Lambda^0$ is a lattice
in $V$ as well.

Throughout this section we keep the above notation.

\begin{lem}
{\rm (\cite[Section 4.1]{Popov})} \label{root} The
following properties hold:
\begin{enumerate}\item[(i)] $\Lambda^0$ is
$G$-invariant; \item[\rm(ii)]
$[\Lambda:\Lambda^0]<\infty$; \item[\rm (iii)]
$\Lambda^0=\Lambda_{1}+\ldots +\Lambda_{n}$;
\item[\rm(iv)] if ${\rm rk}\,(\Lambda) =2n$, then ${\rm
rk}\,(\Lambda^0)=2n$ and ${\rm rk}\,(\Lambda_j)=2$ for
every $j$.
\end{enumerate}
\end{lem}
\begin{proof}
Clearly, $\mathcal L$ is $G$-invariant; whence (i).
Consider the linear operator
\begin{equation}\label{s}
s:=(\I_V-r_1)+\ldots +(\I_V-r_n).
\end{equation}
If $v\in {\rm ker}(s)$, then \eqref{s},
\eqref{formula}, and \eqref{direct} imply $\langle v,
e_1\rangle=\ldots=\langle v, e_n\rangle=0$. Hence
$v=0$, i.e., $s$ is nondegenerate. Therefore ${\rm
rk}\,(\Lambda)={\rm rk}\,(s(\Lambda))$. But \eqref{s},
\eqref{l}, and $G$-invariance of $\Lambda$ imply
$s(\Lambda)\subseteq\Lambda_1+\ldots+\Lambda_n\subseteq
\Lambda^0\subseteq\Lambda$; whence (ii).

 By \cite[Section 3.2]{Popov},
 every
reflection in $G$ is conjugate to a power of some
$r_j$. Hence for every $l\in \mathcal L$, there are
$g\in G$ and an integer $j\in[1, n]$ such that
$g(l)=l_j$. Therefore $g(\Lambda_l)\subseteq
\Lambda_j$. On the other hand, \eqref{l} implies that
$\Lambda_1+\ldots+\Lambda_n$ is invariant with respect
to every $\I_V-r_j$, hence is $G$-invariant. This and
\eqref{G} entails (iii).

By (ii), if ${\rm rk}\,(\Lambda)=2n$, then  ${\rm
rk}\,(\Lambda^0)=2n$. In turn, by (iii) and
\eqref{direct}, the latter equality implies ${\rm
rk}\,(\Lambda_j)=2$ for every $j$. This proves
(iv).\quad $\square$
\renewcommand{\qed}{}\end{proof}

\begin{thm}\label{geom}
Let ${\rm rk}(\Lambda)=2n$. Then:
\begin{enumerate}
\item[\rm(i)] $V/\Lambda^0$ is an abelian variety
isomorphic to a product of mutually isogenous elliptic
curves; \item[\rm(ii)] $V/\Lambda$ is an abelian
variety isogenous to $V/\Lambda^0$, and hence, by
{\rm(i)}, isogenous to a self-product of an elliptic
curve; \item[\rm(iii)] if $G$ is not the
complexification of the Weyl group of an irreducible
root system, then $V/\Lambda$ is isomorphic to a
product of mutually isogenous elliptic curves with
complex multiplication.
\end{enumerate}
\end{thm}
\begin{proof}
Since $\Lambda^0\subseteq\Lambda$,
Lemma~\ref{root}\,(iv)
 implies that $V/\Lambda$,
 $V/\Lambda^0$ are
 complex tori and every
 $l_j/\Lambda_j$ is an elliptic curve.
   From \eqref{direct} and
 Lemma~\ref{root}\,(iii) we deduce that
\begin{equation}\label{product}
V/\Lambda^0\ \text{is isomorphic to}\
l_{1}/\Lambda_{1}\times\ldots
 \times l_{n}/\Lambda_{n},
 \end{equation}
 hence, in particular, $V/\Lambda^0$ is an abelian
 variety. By Lemma~\ref{av}, this implies that
 $V/\Lambda$ is an abelian variety as well and
 \begin{equation}\label{tori}
  \text{abelian varieties}\ V/\Lambda\
  \text{and}\ V/\Lambda^0\ \text{are
  isogenous}.
  \end{equation}

By \eqref{l}, we have $(\I_V-r_k)(l_j)\subseteq l_k$
for every $j$ and $k$. The dimension reason then
implies that the $\C$-linear map
\begin{equation}\label{lin}
(\I_V-r_k)|_{l_j}: l_j\rightarrow l_k
\end{equation}
is either $0$ or an isomorphism. In the latter case,
\eqref{lin} induces an isomorphism of elliptic curves
\begin{equation}\label{iso}
l_j/\Lambda_j\overset{\simeq}\longrightarrow
l_k/\bigl((\I_V-r_k)(\Lambda_j)\bigr).\end{equation} On
the other hand, since $\Lambda$ is $G$-invariant,
$(\I_V-r_k)(\Lambda_j)\subseteq \Lambda_k$. Hence if
\eqref{lin} is an isomorphism, then
$l_k/\bigl((\I_V-r_k)(\Lambda_j)\bigr)$ and
$l_k/\Lambda_k$ are isogenous elliptic curves.  In this
case, by \eqref{iso}, $l_j/\Lambda_j$ and
$l_k/\Lambda_k$ are isogenous elliptic curves as well.
Take now into account that, by \eqref{formula}, the map
\eqref{lin} is an isomorphism if and only if $e_j$ and
$e_k$ are not orthogonal, and, since $G$ is
irreducible, every pair of vectors from the sequence
$e_1,\ldots, e_n$ can be included in a subsequence in
which every two neighboring elements are not
orthogonal. Hence
\begin{equation}\label{jk}
\text{elliptic curves}\ l_j/\Lambda_j\ \text{and}\
l_k/\Lambda_k\ \text{are isogenous for every}\ j, k.
\end{equation}
The proofs of (i) and~(ii) now follow from
\eqref{tori}, \eqref{product}, and \eqref{jk}.

To prove (iii), notice that \eqref{formula} implies
\begin{gather}\label{cyc1}
\begin{gathered}
(\I_V-r_{j_1})(\I_V-r_{j_m}) (\I_V-r_{j_{m-1}})
\ldots(\I_V-r_{j_2})(e_{j_1})=c_{j_1\ldots
j_m}e_{j_1},\hskip 1mm \text{where}\\
\textstyle c_{j_1\ldots j_m}:=\langle e_{j_1},
e_{j_2}\rangle \langle e_{j_2}, e_{j_3}\rangle \ldots
\langle e_{j_{m-1}}, e_{j_m}\rangle \langle e_{j_m},
e_{j_1}\rangle\prod_{t=1}^m(1-\theta_{j_t}).
\end{gathered}
\end{gather}
From \eqref{cyc1} and \eqref{l} we deduce
\begin{equation}\label{tr}
\tr(\I_V-r_{j_1})(\I_V-r_{j_m})(\I_V-r_{j_{m-1}})
\ldots(\I_V-r_{j_2})=c_{j_1\ldots j_m}.
\end{equation}

Additivity of $\tr$ implies that $\Z[\Tr (G)]=\Z[\Tr
(\Z G)]$. Since $\I_V-r_1, \ldots, \I_V-r_n$ generate
the ring $\Z G$, the monomials
$(\I_V-r_{j_1})\ldots(\I_V-r_{j_m})$ generate $\Z G$ as
a $\Z$-module. This and \eqref{tr} entail that $\Z[\Tr
(G)]=\Z[\,\ldots, c_{j_1\ldots j_m},\ldots\,]$,  whence
\begin{equation}\label{=}
\Q(\chi^{}_{G, V})=\Q (\,\ldots, c_{j_1\ldots
j_m},\ldots\,).
\end{equation}

Suppose that $G$
is not the complexification of the Weyl group of an
irreducible root system. Then $\chi^{}_{G, V}$ is not
rational.
Indeed, otherwise Remark~\ref{rgr} and
Theorem~\ref{criterion} would imply 
that $G$ is complexification of a finite real
$n$-dimensional irreducible reflection group that has
an invariant lattice of rank $n$, and it is well known
that such a real group is  the Weyl group of an
irreducible root system, \cite{Bourbaki}.

Since $\chi^{}_{G, V}$ is not rational, \eqref{=}
yields the existence of $j_1,\ldots, j_m$ such that
\begin{equation}\label{neqc}
c_{j_1\ldots j_m}\notin \Q.
\end{equation}
Since $\Lambda$ is $G$-invariant, \eqref{cyc1} implies
that $c_{j_1\ldots j_m}$ is a multiplier of
$\Lambda_{j_1}$, i.e.,
\begin{equation}\label{multiplier}
c_{j_1\ldots j_m}\cdot\Lambda_{j_1}\subseteq
\Lambda_{j_1}.
\end{equation}

Properties \eqref{neqc} and \eqref{multiplier}
imply that $l_{j_1}/\Lambda_{j_1}$ is an elliptic curve
with complex multiplication. But \eqref{tori},
\eqref{product}, and \eqref{jk} imply that $V/\Lambda$
is isogenous to $(l_{j_1}/\Lambda_{j_1})^n$. Thus
$V/\Lambda$ is isogeneous to a self-product of an
elliptic curve with complex multiplication. Now (iii)
follows from Lemma~\ref{Katsura}. \quad $\square$
\renewcommand{\qed}{}\end{proof}

\end{document}